\documentclass[a4paper,12pt]{amsart}
\usepackage{amsthm,vmargin}
\usepackage{amsmath,xypic,amssymb}
\usepackage[colorlinks,pageanchor]{hyperref}
\begin{document}
\newcommand{\Q}{{\mathbb Q}}
\newcommand{\bQ}{\overline{\Q}}
\let\Qb=\bQ
\newcommand{\C}{{\mathbb C}}
\renewcommand{\l}{{\ell}}
\newcommand{\R}{{\mathbb R}}
\newcommand{\Z}{{\mathbb Z}}
\newcommand{\F}{{\mathbb F}}
\newcommand{\gal}{{\rm Gal}}
\newcommand{\N}{{\mathbb N}}
\newcommand{\gQ}{\gal(\bQ/\Q)}
\newcommand{\gl}{{\rm GL}}
\newcommand{\ssl}{{\rm SL}}
\renewcommand{\wp}{{\mathfrak p}}
\renewcommand{\P}{{\mathbb P}}
\renewcommand{\O}{{\mathcal O}}
\newcommand{\Pic}{{\rm Pic\,}}
\newcommand{\Ext}{{\rm Ext}\,}
\newcommand{\rank}{{\rm rk}\,}
\newcommand{\sbull}{{\scriptstyle{\bullet}}}
\newcommand{\bX}{X_{\overline{k}}}
\newcommand{\ch}{\operatorname{CH}}
\newcommand{\tors}{\text{tors}}
\newcommand{\cris}{\text{cris}}
\newcommand{\alg}{\text{alg}}
\newcommand{\tX}{{\tilde{X}}}
\newcommand{\tL}{{\tilde{L}}}
\newcommand{\Hom}{{\rm Hom}}
\newcommand{\spec}{{\rm Spec}}
\let\isom=\simeq
\let\rk=\rank
\let\tensor=\otimes
\newcommand{\X}{\mathfrak{X}}
\newcommand{\mydot}{{\small{\bullet}}}
\let\hom=\Hom
\let\Sl=\ssl
\let\GL=\gl
\let\Gl=\gl

\newtheorem{theorem}[equation]{Theorem}      
\newtheorem{lemma}[equation]{Lemma}          %
\newtheorem{corollary}[equation]{Corollary}  
\newtheorem{proposition}[equation]{Proposition}
\newtheorem{scholium}[equation]{Scholium}

\theoremstyle{definition}
\newtheorem{conj}[equation]{Conjecture}
\newtheorem*{example}{Example}
\newtheorem{question}[equation]{Question}
\theoremstyle{definition}
\newtheorem{remark}[equation]{Remark}

\numberwithin{equation}{section}

\renewcommand{\t}[1]{\tilde{#1}}
\newcommand{\gpb}[1]{(\t{#1},F_i(\t{#1}))}
\title{Infinite Hilbert Class Field Towers from Galois Representations}
\author{Kirti Joshi and Cameron M\lowercase{c}Leman}
\address{Math. department, University of Arizona, 617 N Santa Rita, Tucson, AZ.
85721-0089, USA.} \email{kirti@math.arizona.edu}
\address{Math. department, Willamette University, 800 State Street, Salem, OR.
07304-0089, USA.} \email{cmcleman@willamette.edu}
\date{}


\begin{abstract}
We investigate class field towers of number fields obtained as fixed
fields of modular representations of the absolute Galois group of the
rational numbers.  First, for each $k\in\{12,16,18,20,22,26\}$, we
give explicit rational primes $\l$ such that the fixed field of the
mod-$\l$ representation attached to the unique normalized cusp
eigenforms of weight $k$ on $\Sl_2(\Z)$ has an infinite class field
tower.  Under a conjecture of Hardy and Littlewood, we further prove
that there exist infinitely many such primes for each $k$ (in the above list).  Second,
given a non-CM curve $E/\Q$, we show that there exists an integer $M_E$ such that
the fixed field of the representation attached to the $n$-division
points of $E$ has an infinite class field tower for a set of integers
$n$ of density one among integers coprime to $M_E$.
\end{abstract}
\maketitle

\section{Introduction}
Most current examples of number fields known to have an infinite
(Hilbert) class field tower are constructed ``from the bottom up,''
e.g., by beginning with a fixed number field and constructing
extensions of that number field in which a large number of primes
ramify.  If sufficiently many primes ramify, invoking results from
genus theory and one of many variants of the Golod-Shafarevich
inequality is enough to prove the class field tower infinite.  We
investigate instead class field towers of number fields obtained
``from the top down,'' i.e., as fixed fields of representations of the
absolute Galois group $\gal(\Qb/\Q)$ and prove that many naturally
occurring fields of arithmetic interest have infinite Hilbert class
field towers. In the process, we prove also prove, assuming a
conjecture of Hardy-Littlewood, that there are infinitely cyclotomic
fields of prime conductor with infinite Hilbert class field towers.

We consider two primary sources for these representations.  First, 
we consider Galois representations arising from modular forms.  For
any $k\in\{12,16,18,20,22,26\}$, it is well-known that there exists a
unique normalized cuspidal Hecke eigenform $\Delta_k$ on $\Sl_2(\Z)$
of weight $k$ with integer Fourier coefficients. 
For each such form and rational prime $\l$, we consider the fixed field of the associated
mod-$\l$ representation $\rho_{\Delta_k,\l}$.  
We find explicit examples of primes $\l$ such that each of these fixed
fields have an infinite class field tower.  The key new idea is to use
certain auxiliary cubic fields introduced by Daniel Shanks (see
\cite{shanks74}). Using these auxiliary fields and a refined
Golod-Shafarevich type inequality of R.~Schoof (see \cite{schoof86})
we show that, for suitable primes $\ell$, the fixed fields of the
mod-$\ell$ representations alluded to earlier have an infinite Hilbert
class field tower.  Moreover, assuming a well-known conjecture of
Hardy and Littlewood on prime values of quadratic polynomials, we
prove the existence of infinitely many such $\ell$ (see
Theorem~\ref{cyc} and Corollary~\ref{Cor}).  Further, the fields
arising from these mod-$\l$ representations have Galois groups
containing $\Sl_2(\Z/\l)$ and are ramified at a single finite prime
(in contrast to the number fields shown to have an infinite class
field tower via genus theory). As a consequence, we also show that the
conjecture of Hardy and Littlewood implies the existence of infinitely
many primes $\ell$ such that the cyclotomic fields $\Q(\zeta_\ell)$
have infinite Hilbert class field towers.

For the second construction, let $E/\Q$ be an elliptic curve without
complex multiplication.  For any $n\geq 1$, the absolute Galois group
$\gal(\Qb/\Q)$ acts on the $n$-torsion points
$E[n](\Qb)\cong(\Z/n\Z)^{2}$ of $E$.  The fixed field $K_n$ of the
kernel of the associated representation
$\rho_n:\gal(\Qb/\Q)\rightarrow\gl_2(\Z/n)$ is of significant
arithmetic interest.  In Theorem~\ref{SmallThm}, we give an elementary
proof that for almost all integers $n$ coprime to a fixed integer
$A_E$ (almost all here means outside a set of density zero), the field
$K_n$ has an infinite class field tower.  Further, known information
about these fixed fields can be translated into information about the
fields arising in these towers.  Finally, we note that the
construction can be made fairly explicit--in Remark \ref{FurutaNote},
we give a specific example of an elliptic curve and an integer $n$
provided by the result.

We would like to thank William McCallum and Dinesh Thakur for comments.

\section{Infinite Class Field Towers from Hecke Eigenforms}
For each even integer $k$ with $12\leq k\leq 26$ and $k\neq 14,24$,
there exists a unique normalized cuspidal Hecke eigenform on
$\ssl_2(\Z)$ of weight $k$, which we denote here by $\Delta_k$.  The
theorem of Deligne-Serre provides for each $\Delta_k$ and prime $\ell$
a Galois representation $\rho_{\Delta_k,\ell}:\gal(\Qb/\Q)\to
\gl_2(\mathbb{F}_{\ell})$ which is unramified outside $\ell$. This
representation is the reduction modulo $\ell$ of the $\ell$-adic
representation associated to $\Delta_k$.  Let $K_{\Delta_k,\ell}$ be
the fixed field of the representation.  By a well-known result of
Serre and Swinnerton-Dyer (see \cite{serre73}), for any $\Delta_k$ and
sufficiently large prime $\ell$, the Galois group
$\gal(K_{{\Delta_k},\ell}/\Q)$ is the subgroup
$$G_{\ell}=\left\{g\in \gl_2(\mathbb{F}_{\ell}):\det(g)\in
(\mathbb{F}_{\ell}^*)^{k-1}\right\}$$ of $\gl_2(\mathbb{F}_{\ell})$.
Note that if $(k-1,\ell-1)=1$, we obtain
$G_{\ell}=\gl_2(\mathbb{F}_{\ell})$.  The theorem is effective and one
has a complete list of the exceptional primes:
\begin{center}
\begin{tabular}{|c|c|}\hline
$k$&exceptional $\ell$ for $\Delta_k$\\\hline
$12$&$2,3,5,7,23,691$\\\hline
$16$&$2,3,5,7,11,31,59,3617$\\\hline
$18$&$2,3,5,7,11,13,43867$\\\hline
$20$&$2,3,5,7,11,13,,283,617$\\\hline
$22$&$2,3,5,7,13,17,131,593$\\\hline
$26$&$2,3,5,7,11,17,19,657931$\\\hline
\end{tabular}
\end{center}
The field $K_{\Delta_k,\ell}=\Qb^{\rho_{\Delta_k,\ell}}$ certainly
contains the fixed field $\Qb^{\det(\rho_{\Delta_k,\ell})}$ of the
determinant representation of $\rho_{\Delta_k,\ell}$. This is a field
which is unramified outside $\ell$ and has
$(\mathbb{F}_{\ell}^*)^{k-1}$ as its Galois group over $\Q$. If
$(k-1,\ell-1)=1$, then the Galois group has order $\ell-1$ and the
extension is unramified outside $\ell$. By class field theory, such an
extension corresponds to a suitable quotient of $\Z_\ell^*$. The only
quotient of this group of order $\ell-1$ is the quotient by the
subgroup of $\Z_\ell^*$ consisting of $\ell$-adic units congruent to
$1$ mod $\ell$, and thus the determinantal fixed field is
$\Q(\zeta_\ell)$.  Since any extension of a field with an infinite
class field tower itself has an infinite class field tower (a variant
of this easy lemma is proved in the middle of the proof of Theorem
\ref{SmallThm}), we have proven the following:
\begin{theorem}\label{BigThm}
Let $k\in\{12,16,18,20,22,26\}$, and let $\ell$ be a prime such that:
\begin{itemize}
\item $\l$ is not exceptional for $\Delta_k;$
\item $(k-1,\l-1)=1;$
\item $\Q(\zeta_\l)$ has an infinite class field tower.
\end{itemize}
Then $\gal(K_{\Delta_k,\ell}/\Q)=\gl_2(\mathbb{F}_{\ell})$, and
$K_{\Delta_k,\l}$ has an infinite class field tower.
\end{theorem}

We now turn to searching for primes which satisfy the hypotheses of
the theorem.  As a first example, we note that the prime $\ell=877$
satisfies the conditions of the theorem for the Ramanujan form
$\Delta_{12}$.  Namely, we have that
$\operatorname{gcd}(12-1,877-1)=1$ and that $\ell=877$ is not an
exceptional prime for $\Delta_{12}$, so $\rho_{\Delta_{12},877}$ is
surjective.  Moreover, it is shown in \cite{schoof86} that
$\Q(\zeta_{877})$ has an infinite Hilbert class field tower.
Theorem~\ref{BigThm} then gives that $K_{\Delta_{12},877}$ has an
infinite class field tower.

What principally remains to do is to generalize Schoof's argument to
produce a large class of primes $\l$ such that $\Q(\zeta_\l)$ has an
infinite class field tower.  We do this below, and show that if we
assume the conjecture of Hardy and Littlewood presented below, the set
of such primes is in fact infinite.  The conjecture arose from their
famous ``circle method,'' and should be viewed as the quadratic
analogue of Dirichlet's Theorem on primes in arithmetic (i.e.,
``linear'') progressions.
\begin{conj}[Hardy-Littlewood, \cite{HL23}]\label{HL}
If $h(x):=ax^2+bx+c\in\Z[x]$ satisfies:
\begin{itemize}
\item the quantities $a+b$ and $c$ are not both even;
\item the discriminant $D(h):=b^2-4ac$ is not a square;
\end{itemize}
then $h(x)$ represents infinitely many prime values.
\end{conj}
We now show that the conjecture implies the existence of infinitely
many cyclotomic fields \emph{of prime conductor} with an infinite
class field tower.  In contrast, most techniques used to construct
infinite class field towers (note in particular Remark
\ref{FurutaNote}) provide fields with highly composite conductors.
\begin{theorem}\label{cyc}
Under the assumption of Conjecture \ref{HL}, there exist infinitely
many primes $\l$ such that $\Q(\zeta_\l)$ has an infinite class field
tower.
\end{theorem}
\begin{proof}
Let $k$ be a positive integer and let $m=12k+2$.  Consider the cubic
polynomial
$$f_m(x)=x^3-mx^2-(m+3)x-1,$$ with discriminant $(m^2+3m+9)^2$.  With
notation as in the Hardy-Littlewood conjecture, the quadratic
polynomial $m^2+3m+9=144k^2+84k+19$ has (viewed as a polynomial in
$k$) odd $c$ and the non-square discriminant
$$D(144k^2+84k+19)=-3888.$$

Thus, assuming the conjecture, there are infinitely many values of $k$
for which $m^2+3m+9=:\l$ is prime.  For the remainder of the proof, we
restrict to such $k$, $m$, and $\l$.  In this case, the splitting
field $F_m$ of $f_m(x)$ is one of Shanks' ``simplest cubic fields,'' a
totally real cyclic cubic field of prime conductor $\l$ (see
\cite{shanks74}).  Thus $F_m$ is the unique cubic subfield of
$\Q(\zeta_\l)$.  Now note that since $m\equiv 2$ mod $12$, we have
$\l\equiv 7$ mod $12$.  Since $\l\equiv 1$ mod $6$, there is a unique
sextic subfield of $\Q(\zeta_\l)$, which we denote by $L_m$. Further,
$6\nmid \frac{\l-1}{2}$, so $L_m\not\subset\Q(\zeta_\l)^+$ is totally
imaginary.  Let $h$ be the class number of $F_m$ and
$L_m'=L_mF_m^{(1)}$, giving the following diagram of fields:
\begin{equation*}
\xymatrix@R=1pc@C=3pc{
  \Q(\zeta_\l)\ar@{-}[dd]_{2}&&L_m'\ar@{-}[dd]^2\\
  &L_m\ar@{-}[ul]_(.5){\frac{\l-1}{6}}\ar@{-}[dd]_2\ar@{-}[ur]^(.5){h}\\
\Q(\zeta_\l)^+\ar@{-}[dd]_{\frac{\l-1}{2}}&&F_m^{(1)}\\
&F_m\ar@{-}[ul]\ar@{-}[dl]^3\ar@{-}[ur]_{h}\\
\Q  &}
\end{equation*}
Denote by $d_2E_K$ the 2-rank of the unit group of a number field
$K$. Applying \cite[Proposition 3.3]{schoof86} to the extension
$L_m'/F_m^{(1)}$, a sufficient condition for $L_m'$ to have an
infinite class field tower is that that the number $\rho$ of ramified
primes in this extension satisfies
$$\rho\geq 3+d_2E_{F_m^{(1)}}+2\sqrt{d_2E_{L_m'}+1}.$$

Since $L_m'$ is totally imaginary and $F_m^{(1)}$ is totally real,
Dirichlet's Unit theorem easily calculates the right-hand side of this
inequality to be $3+3h+2\sqrt{3h+1}$.  Now we count ramified primes:
First, all $3h$ infinite places of $F_m^{(1)}$ ramify in
$L_m'$. Second, since $L_m/F_m$ is totally ramified at the prime
$\lambda$ of $F_m$ above $\l$, and $\lambda$ (being principal) splits
completely in $F_m^{(1)}$, the extension $L_m'/F_m^{(1)}$ is ramified
also at the $h$ primes of $F_m^{(1)}$ above $\l$.  In sum, this gives
$\rho=4h$, and we find that the equality is satisfied whenever $h\geq
18$.  Using that $h\rightarrow\infty$ as $m\rightarrow \infty$ (see
\cite{shanks74}), we now get infinitely many $m$ such that $L_m'$ has
an infinite class field tower.  Note that $F_m^{(1)}$ and $L_m$ are
linearly disjoint over $F_m$ since the former is unramified, and the
latter is totally ramified at primes above $\l$.  By disjointness,
$L_m'/L_m$ is abelian and unramified, and so $L_m$, and consequently
$\Q(\zeta_\l)$, also have infinite class field towers.
\end{proof}

\begin{corollary}\label{Cor}
Assume Conjecture $\ref{HL}$. Then for each
$k\in\{12,16,18,20,22,26\}$, there are infinitely many primes $\l$
such that $K_{\Delta_k,\l}$ has an infinite class field tower.
\end{corollary}
\begin{proof}
It is easy to verify that primes $\l$ of the form $m^2+3m+9$ are never
congruent to 1 mod $(k-1)$ for each $k$ in the list, so $(\l-1,k-1)=1$
for any $\l$ constructed by the theorem.  Avoiding the finitely many
primes in the table given before Theorem~\ref{BigThm} for which the
representation is not surjective, the remaining primes constructed in
Theorem~\ref{cyc} satisfy all of the hypotheses of
Theorem~\ref{BigThm}.
\end{proof}

\begin{remark}
Hardy and Littlewood also provide an asymptotic version of their
conjecture (see again \cite{HL23}).  Applied to the polynomial $h(k)$
used in the proof, the number $P_h(x)$ of primes less than $x$ which
are represented by $h(k)=144k^2+84k+19$ satisfies
$$
P_h(x)\sim \frac{1}{4}\prod_{p=5}^\infty\left(1-\frac{\left(\frac{-3888}{p}\right)}{p-1}\right)\frac{\sqrt{x}}{\log x}\approx \frac{.28\sqrt{x}}{\log x},
$$ where we have used SAGE to approximate the constant by including
the terms in the product for all primes $p\leq 10^7$.  This thus also
provides an asymptotic lower bound for the number of primes $\l\leq x$
such that $\Q(\zeta_\l)$ has an infinite class field tower.
\end{remark}

Finally, we remark that setting $m=12k+2$ in the proof of Theorem
\ref{BigThm} was overly restrictive on our choice of $m$, designed
only to ensure that $\l\equiv 7\mod 12$.  This can be equally well
achieved by insisting that $m\equiv 2,7,10,\text{ or }11\mod 12$.
Searching Shanks' Table 1 for such values of $m$ giving $h\geq 18$
provides the first few examples of $\Q(\zeta_\l)$ provided by the
proof: $\l\in\{2659,3547,5119,8563,\ldots\}.$ Note that for each of
these specific values of $\l$, Theorem \ref{BigThm} is proved
unconditionally, the Hardy-Littlewood conjecture being used only to
guarantee that there are infinitely many primes in this list.

\section{Infinite Class Field Towers from non-CM Elliptic Curves}
A second class of fields of arithmetic interest we discuss are the
fixed fields of representations attached to elliptic curves. Let $E$
be an elliptic curve without complex multiplication.  Let
$\rho_n:\gQ\to\gl_2(\Z/n)$ be the Galois representation associated to
the $n$-torsion points of $E$.  Let $A_E$ be the product of all the
exceptional primes of $E$, the finite set of primes $\l$ such that
$\rho_\l$ is not surjective.  Then for all $n$ relatively prime to
$M_E:=30A_E$, the representation $\rho_n$ is surjective, and hence
$\gal(K_n/\Q)\cong\gl_2(\Z/n)$ (see
\cite{serre72,serre-abelian,kani05}).  We note that this set of
primes, and hence the constants $A_E$ and $M_E$, have been studied
extensively.  In particular, we note that explicit bounds on these
constants are known \cite[Theorem A.1 and Theorem 2]{Co05}), and
$M_E=30$ for almost all elliptic curves (see Remark~\ref{DukeNote}).
For a number field $F$, let $F^{(m)}$ denote the $m$-step in the
Hilbert class field tower over $F$.

\begin{theorem}\label{SmallThm}
Let $E/\Q$ be an elliptic curve without complex multiplication, and
let $\rho_n$, $K_n$, $A_E$, and $M_E$ be as in the preceding
paragraph.  Let $S$ be the set of integers prime to $M_E$.  Then for
all $n\in S$ outside a subset of density zero, the field $K_n$ has an
infinite Hilbert class field tower.  Furthermore, for such $n$, there
is a natural surjection of class groups
$\operatorname{Cl}(K_n^{(m)})\rightarrow
\operatorname{Cl}(\Q(\zeta_n)^{(m)})$ for each $m\geq 1$.

\end{theorem}
\begin{proof}
For $n$ prime to $M_E$, the representation $\rho_n$ is surjective, and so
$\text{Image}(\rho_n)\cong\gal(K_{n}/\Q)\cong\gl_2(\Z/n).$ The fixed
field of the kernel of the composite representation $\gQ\to(\Z/n)^*$
arising from the exact sequence
$$ 1\to\ssl_2(\Z/n)\to\gl_2(\Z/n)\to(\Z/n)^*\to 1
$$ is the $n$-th cyclotomic field $\Q(\zeta_n)$, and so we have an
inclusion of fields $\Q\subset \Q(\zeta_n)\subset K_n$. By a result of
Shparlinski (see \cite{shparlinski08}), the set of $n$ coprime to
$M_E$ for which $\Q(\zeta_n)$ has an infinite class field tower has
density one in the set of integers coprime to $M_E$. For the remainder
of the proof, $n$ will denote such an integer.  We claim that the
fields $\Q(\zeta_n)^{(m)}$ and $K_{n}$ are linearly disjoint
extensions of $\Q(\zeta_n)$ for all $m\geq 1$. Let $H_{n,m}$ denote
their intersection.  Consider the lattice diagram of fields:
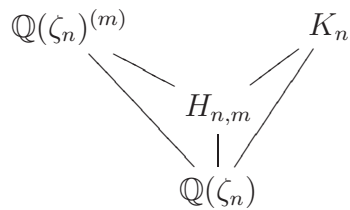
\begin{equation*}
\xymatrix@R=1pc@C=1pc{
  \text{   }\Q(\zeta_n)^{(m)}\ar@{-}[ddr] \ar@{-}[rd]_{}  & & K_n\ar@{-}[ld]\ar@{-}[ldd]  \\
  & H_{n,m}\ar@{-}[d] &\\
   & \Q(\zeta_n)  &}
\end{equation*}
Since $\Q(\zeta_n)^{(m)}/\Q(\zeta_n)$ is a Galois extension with
solvable Galois group (being constructed via a series of abelian
extensions), so is $H_{n,m}/\Q(\zeta_n)$.  But $\gal(K_n/\Q)\cong
\ssl_2(\Z/n)$ is perfect for $(n,30)=1$ (see \cite{serre-abelian}), and
so has no non-trivial solvable quotients.  Thus
$\gal(H_{n,m}/\Q(\zeta_n))$ is trivial and the two fields are linearly
disjoint.  Let $L_{n,m}$ be the compositum $K_n\Q(\zeta_n)^{(m)}$,
giving the following diagram for each $m\geq 1$.
\begin{equation*}
  \xymatrix@R=.5pc{
  L_{n,m}\ar@{-}[dd]\ar@{-}[dr]&\\ &K_n\ar@{-}[dd]\\ \Q(\zeta_n)^{(m)}\ar@{-}[dr]&
  \\ &\Q(\zeta_n)}
\end{equation*}
By linear disjointness, we have $\gal(L_{n,m}/\Q(\zeta_n)^{(m)})\cong
\gal(K_n/\Q(\zeta_n))\cong \ssl_2(\Z/n)$ for all $m\geq 0$.  Further,
each extension $L_{n,m}/L_{n,m-1}$ is unramified abelian.  Thus for
each $m\geq 0$, $L_{n,m}$ is an unramified solvable extension of
$L_{n,0}=K_n$, and hence is contained in the Hilbert class field tower
over $K_n$.  Since the degrees of $L_{n,m}$ go to infinity with $m$,
we see that $K_n$ has an infinite class field tower, as desired.
Finally, $L_{m,n}$ is an unramified extension of $K_n$ with
$\gal(L_{m,n}/K_n)$ solvable of derived length $\leq m$, so is
contained in the maximal such extension $K_n^{(m)}$.  The restriction
map on Galois groups
\begin{equation*}
  \xymatrix@R=.5pc{
\gal(K_n^{(m+1)}/K_n^{(m)})\ar@{->>}[r]&\gal(K_n^{(m)}\Q(\zeta_n)^{(m)}/K_n^{(m)})\cong \gal(\Q(\zeta_n)^{(m+1)}/\Q(\zeta_n)^{(m)})
}
\end{equation*}
corresponds by class field theory to the desired surjection
$\operatorname{Cl}(K_n^{(m)})\rightarrow
\operatorname{Cl}(\Q(\zeta_n)^{(m)})$ of class groups.
\end{proof}
\begin{remark}\label{FurutaNote}
The density result of Shparlinski used in the proof is based on an
explicit construction due to Furuta (see \cite{furuta72}, Theorem
4). Namely, we can find explicit examples of the $n$ described in the
theorem by choosing a rational prime $\l$ and taking $n$ to be a
product of nine or more rational primes congruent to 1 mod $\l$ and
prime to $M_E$.  The field $K_n$ then has an infinite $\l$-class
field tower. For example, consider the elliptic curve
$$E:y^2+y=x^3-x$$ of conductor $37$.  By \cite[5.5.6, Page
  310]{serre72} we find that $A_E=1$ and so $M_E=30$. We take $n$ to
be the product of the first nine primes which are congruent to $1$
mod $5$:
$$n=11\cdot31\cdot41\cdot61\cdot71\cdot101\cdot131\cdot151\cdot181.$$
Then $K_{n}$ has an infinite Hilbert $5$-class field tower and is
unramified outside primes dividing $37n$.
\end{remark}

\begin{remark}\label{DukeNote}
The constant $A_E$ in Serre's Theorem has been studied extensively.
In \cite{duke97} it was shown that almost all elliptic curves over
$\Q$ have Serre constant $A_E=1$ (and thus $M_E=30$). Thus for almost
all elliptic curves over $\Q$ and almost all integers $n$ with
$(n,30)=1$, the field $K_n$ has an infinite class field tower.
\end{remark}

Finally, we note that variants of these techniques apply to other
arithmetically significant fields.  For example, let $E$ be a
semistable elliptic curve over $\Q$ of prime conductor $\l\geq11$, and
suppose that $\Q(\zeta_\l)$ has an infinite Hilbert class field tower.
By the above techniques, $K_{E,\l}$ has an infinite Hilbert class
field tower.  Further, $K_{E,\l}$ is unramified outside $\l$, and and
by a well-known theorem of Mazur (see \cite{mazur78}), we have
$\gal(K_{E,\l}/\Q)=\Gl_2(\Z/\l)$.  Unfortunately, it is unknown if
there exists an elliptic curve of conductor $\l$ for infinitely many
primes $\l$ (though this is widely believed). Regardless, an
examination of Cremona's tables of elliptic curves (see
\cite{cremona-tables}) shows that there are no elliptic curves with
conductor $877$ but that there exist elliptic curves with conductor
$3547$, the second prime provided in the discussion after Corollary
\ref{Cor}. Thus we can find examples of elliptic curves to which this
variant approach applies.

\bibliographystyle{amsplain}

\begin{thebibliography}{99}
\bibitem{Co05} Cojocaru, Alina Carmen.  \emph{On the surjectivity of
  Galois representations associated to non-CM elliptic curves.}  Canad. Math. Bull., {\bf 48}(1):16--31, 2005.
\bibitem{cremona-tables} John~Cremona,  Algorithms for modular elliptic curves
\href{http://www.warwick.ac.uk/staff/J.E.Cremona/book/fulltext/index.html}{online
edition}.
\bibitem{duke97}W. D. Duke, \emph{Elliptic curves with no exceptional primes}, C. R. Math. Acad.
Sci. Paris S�er. I 325 (1997), pp. 813�-818.
\bibitem{furuta72} Yoshiomi Furuta.  \emph{On Class Field Towers and
  the Rank of Ideal Class Groups}.  Nagoya Math J, {\bf 48}, 1972.
  147-157.
\bibitem{HL23} G.H. Hardy and J. E. Littlewood.  \emph{Some problems
  of `Partitio numerorum'; III: On the expression of a number as a sum
  of primes.}  Acta Math. 44 (1923), no. 1, 1--70.
\bibitem{kani05} Kani, Ernst.  \emph{Appendix to: On the surjecticity of
  Galois representations associated to non-CM elliptic curves.}  Canad. Math. Bull., {\bf 48}(1):16--31, 2005.
  \bibitem{mazur78}
  Barry Mazur. \emph{Rational Isogenies of prime degree}, Invent. Math. 44 (1978),
pp. 129--162.
\bibitem{schoof86} Rene Schoof, \emph{Infinite class field towers of
  quadratic fields}, J. Reine Angew. Math. 372 (1986), 209--220.
\bibitem{serre-abelian} Jean-Pierre Serre, \emph{Abelian $\ell$-adic representations and elliptic curves}, McGill University lecture
   notes written with the collaboration of Willem Kuyk and John Labute, W. A. Benjamin, Inc.,
   New York-Amsterdam, 1968.
\bibitem{serre72} ---, \emph{Propri\'{e}t\'{e}s Galoisiennes des
  points d'ordre fini des courbes elliptiques}, Inventiones
  Mathematicae {\bf 15}, 1972, pp. 259-331.
\bibitem{serre73} ---, \emph{Congruences et formes modulaires}
  [d'apr\`{e}s H. P. F. Swinnerton-Dyer], S\'{e}minaire Bourbaki 24e
  ann\'{e}e (1971/1972), Exp. No. 416, Springer, Berlin, 1973,
  pp. 319–338. Lecture Notes in Math., Vol. 317.
\bibitem{shanks74} Daniel Shanks.  \emph{The Simplest Cubic Fields.}
  Mathematics of Computation, Volume 28, Number 128.  1974.
  1137-1152.
\bibitem{shparlinski08} Igor E. Shparlinski, \emph{Infinite Hilbert class field towers
  over cyclotomic fields}, Glasg. Math. J. 50 (2008), no. 1, 27–32.
\end{thebibliography}

\end{document}